\documentclass[11pt]{article}
\usepackage{amsfonts,amsmath,amsthm,epsfig,verbatim,amssymb,amscd,eucal,bbm}

\usepackage[latin1]{inputenc}

\theoremstyle{plain}
\newtheorem{thm}{Theorem}[section]

\newtheorem{lem}[thm]{Lemma}

\newtheorem{defn-prop}[thm]{Definition-Proposition}

\newcommand{\cP}{{\cal P}}

\newcommand{\N}{\mathbb{N}}

\newcommand{\dist}{{r}}

\begin{document}

\noindent

\title{Efficient routing in Poisson small-world networks}

\author{ M. {\sc Draief} \thanks{Statistical Laboratory, Centre for
Mathematical Sciences, Wilberforce Road, Cambridge CB3 0WB UK
E-mail: {\tt M.Draief@statslab.cam.ac.uk}} ~and~A. {\sc Ganesh}
\thanks{Microsoft Research, 7 J.J. Thomson Avenue, Cambridge CB3 0FB E-mail: {\tt ajg@micorsoft.com} }}

\date{}

\maketitle

\begin{abstract}
\noindent

In recent work, Jon Kleinberg considered a small-world network
model consisting of a $d$-dimensional lattice augmented with
shortcuts. The probability of a shortcut being present between two
points decays as a power, $r^{-\alpha}$, of the distance $r$
between them. Kleinberg showed that greedy routing is efficient if
$\alpha = d$ and that there is no efficient decentralised routing
algorithm if $\alpha \neq d$. The results were extended to a
continuum model by Franceschetti and Meester. In our work, we
extend the result to more realistic models constructed from a
Poisson point process, wherein each point is connected to all its
neighbours within some fixed radius, as well as possessing random
shortcuts to more distant nodes as described above.
\end{abstract}

%%% ----------------------------------------------------------------------

\section{Introduction}

A classical random graph model introduced by Erd\H{o}s and R\'enyi
consists of $n$ nodes, with the edge between any pair of vertices
being present with probability $p(n)$, independent of other pairs.
Recently, there has been  considerable interest in alternative
models where the nodes are given coordinates in an Euclidean
space, and the probability of an edge between a pair of nodes
$u$ and $v$ is given by a function $g(\cdot)$ of the distance $r(u,v)$
between the nodes; edges between different node pairs are independent.
Such `random connection' or `spatial random graph' models and variants
thereof arise, for instance, in the study of wireless communication networks.

\medskip

The ``small-world phenomenon" (the principle that all people are linked
by short chains of acquaintances), which has long been a matter of
folklore,  was inaugurated as an area of experimental study in the
social sciences through the pioneering work of Stanley Milgram \cite{Mil67}.
Recent works have suggested that the phenomenon is pervasive in networks
arising in nature and technology, and motivated interest in mathematical
models of such networks. While Erd\H{o}s-R\'enyi random graphs possess
the property of having a small diameter (smaller than logarithmic in
the number of nodes, above the connectivity threshold for $p(n)$),
they are not good models for social networks because of the independence
assumption. On the other hand, spatial random graphs are better at
capturing clustering because of the implicit dependence between edges
induced by the connection function $g(\cdot)$.

\medskip

Watts and Strogatz \cite{WaSt98} conducted a set of re-wiring experiments
on graphs, and observed that by re-wiring a few random links in finite
lattices, the average path length was reduced drastically (approaching
that of random graphs). This led them to propose a model of ``small-world
graphs" which essentially consists of a lattice augmented with random
links acting as shortcuts, which play an important role in shrinking
the average path link. By the length of a path we mean the number of
edges on it, and distance refers to graph distance (length of shortest path)
unless otherwise specified.

\medskip

The diameter of the Watts-Strogatz model in the 1-dimensional case was
obtained by Barbour and Reinert \cite{BR}. Benjamini and Berger \cite{BB}
considered a variant of this 1-dimensional model wherein the shortcut
between any pair of nodes, instead of being present with constant
probability, is present with probability given by a connection function
$g(\cdot)$; they specifically considered connection functions of the form
$g(r) \sim \beta r^{-\alpha}$, where $\beta$ and $\alpha$ are given constants,
and $r(u,v)$ is the graph distance between $u$ and $v$ in the underlying
lattice (i.e., the $L_1$ distance).

\medskip

The general $d$-dimensional version of this model, on the finite lattice
with $n^d$ points, was studied by Coppersmith et al. \cite{CGS}. They
showed that the diameter of the graph is (i) $\Theta(\log n/\log \log n)$
if $\alpha=d$, (ii) at most polylogarithmic in $n$ if $d<\alpha<2d$, and
(iii) at least polynomial in $n$ if $\alpha>2d$. Finally, it was shown by
Benjamini et al.  \cite{BKPS} that the diameter is a constant if $\alpha<d$.
%Note that in this model, the mean degree
%of a node is $\Theta(1)$ if $\alpha>d$, $\Theta(\log n)$ if $\alpha=d$ and
%$\Theta(n^{\alpha})$ for some $\alpha>0$ if $\alpha<d$. This raises the
%question of whether the decrease in diameter is due to the greater
%number of shortcuts or better properties of the shortcuts.

\medskip

The sociological experiments of Milgram demonstrated not only that
there is a short chain of acquaintances between strangers but also
that they are able to find such chains. What sort of graph models
have this property? Specifically, when can decentralised routing
algorithms (which we define later) find a short path between
arbitrary source and destination nodes?

\medskip

This question was addressed by Jon Kleinberg \cite{Klei00} for the
class of finite $d$-dimensional lattices augmented with shortcuts,
where the probability of a shortcut being present between two nodes
decays as a power, $r^{-\alpha}$ of the distance $r$ between them.
Kleinberg showed that greedy routing is efficient if $\alpha = d$
and that there is no efficient decentralised routing algorithm if
$\alpha \neq d$. The results were extended to a continuum model by
Franceschetti and Meester \cite{FrMee04}. Note that these results
show that decentralised algorithms cannot find short routes when
$\alpha \neq d$, even though such routes are present for $\alpha<2d$
by the results of Benjamini et al. and Coppersmith et al. cited
above; when $\alpha > 2d$, no short routes are present.

\section{Our Model}

In this work, we consider a model constructed from a Poisson point
process on a finite square, wherein each point is connected to all
its neighbours within some fixed radius, as well as possessing random
shortcuts to more distant nodes. More precisely:
\begin{itemize}
\item We consider a sequence of graphs indexed by $n\in \N$.
\item Nodes form a Poisson process of rate $1$ on the square
$[0,\sqrt{n}]^2$.
\item Each node $x$ is linked to all nodes that are distance less
that $r_n=\sqrt{c\log n}$ for a sufficiently large constant, $c$.
In particular, if $c>1/\pi$, then this graph is connected with high
probability (abbreviated whp, and meaning with probability going to
1 as $n$ tends to infinity); see \cite{penrose}.
These links are referred to as local edges and the corresponding nodes
as the local contacts of $x$.
\item For two nodes $u$ and $v$ such that $\dist(u,v) > \sqrt{c\log n}$,
the edge $(u,v)$ is present with probability $a_n\dist(u,v)^{-\alpha}
\wedge 1$. Such edges are referred to as shortcuts. The parameter $a_n$
is chosen so that the expected number of shortcuts per node is equal to
some specified constant, ${\overline d}$.
\end{itemize}

The objective is to route a message from an arbitrary source node $s$
to an arbitrary destination $t$ using a small number of hops. We are
interested in decentralised routing algorithms, which do not require
global knowledge of the graph topology. It is assumed throughout that
each node knows its location (co-ordinates) on the plane, as well as
the location of all its neighbours, both local and shortcut, and of
the destination $t$. We show that efficient decentralised routing is
possible only if $\alpha = 2$. More precisely, we show the following:
\begin{itemize}
\item[$\bullet$] $\alpha=2$: there is a greedy decentralised algorithm
to route a message from source to destination in $O(\log^2 n)$ hops.
\item[$\bullet$] $\alpha<2$: any decentralised routing needs more than
$n^{\gamma}$ hops on average, for any $\gamma$ such that
$\gamma<(2-\alpha)/6$.
\item[$\bullet$] $\alpha>2$: any decentralised routing needs more than
$n^{\gamma}$ hops on average, for any $\gamma<\frac{\alpha-2}{2(\alpha-1)}$.
\end{itemize}

As noted by Kleinberg for the lattice model, the case $\alpha=2$
corresponds to a ``scale-free" network: the expected number of shortcuts
from a node $x$ to nodes which lie between distance $r$ and $2r$ from it
is the same for any $r$. It was observed by Franceschetti and Meester in
their continuum model that this property is related to the impossibility
of efficient decentralised routing when $\alpha \neq 2$ through the fact
that shortcuts can't make sufficient progress towards the destination
when $\alpha > 2$ (they are too short) while they can't home in on small
enough neighbourhoods of the destination when $\alpha < 2$ (they are too
long). Similar remarks apply to our model as well.

\medskip

A model very similar to ours was considered by Sharma and Mazumdar
\cite{SM} who use it to describe an ad-hoc sensor network. The
sensors are located at the points of a Poisson process and can
communicate with nearby sensors through wireless links (corresponding
to local contacts). In addition, it is possible to deploy a small number
of wired links (corresponding to shortcuts), and the question they
address is that of how to place these wired links in order to enable
efficient decentralised routing.

\medskip

In the analysis presented below, we ignore edge effects for ease of
exposition. This is equivalent to considering distances as being defined
on the torus obtained by identifying opposite edges of the square.

\section{Efficiency of greedy routing when $\alpha=2$}

When $\alpha=2$, we show that the following \emph{approximately}
greedy algorithm succeeds whp in reaching the destination in a number
of hops which is polylogarithmic in $n$, the expected number of nodes.

\medskip

Denote by $C(u,r)$ the circle of radius $r$ centred at node $u$.
If there is no direct link from the source $s$ to the destination
$t$, then the message is passed via intermediate nodes as follows.
At each stage, the message carries the address (co-ordinates) of
the destination $t$, as well as a radius $r$ which is initialised
to $\dist(s,t)$, the distance between $s$ and $t$. Suppose the message
is currently at node $x$ and has radius $r > \sqrt{c\log n}$.
(If $r \le \sqrt{c\log n}$, then the node which updated $r$ would have
contained $t$ in its local contact list and delivered the message
immediately.) If node $x$ has a shortcut to some node $y \in A(t,r)$,
where the annulus $A(t,r)$ is defined as $A(t,r)= C(t,\frac{r}{2})
\setminus C(t,\frac{r}{4})$, then $x$ forwards the message to
$y$. If there is more than one such node, the choice can be arbitrary.
Otherwise, it forwards the message to one of its local
contacts which is closer to $t$ than itself. When a node $y$
receives a message, it updates $r$ to $r/2$ if
$\dist(y,t)\le r/2$, and leaves $r$ unchanged otherwise.

\medskip

In other words, if $x$ can find a shortcut which reduces the distance to
the destination by at least a half but by no more than three-quarters, it
uses such a shortcut. Otherwise, it uses a local contact to reduce the
distance to the destination. In that sense, the algorithm is approximately
greedy. The reason for considering such an algorithm rather than a greedy
algorithm that would minimize the distance to the destination at each
step is to preserve independence, which greatly simplifies the analysis.
Note that if a greedy step from $x$ takes us to $y$ (i.e., of all nodes
to which $x$ possesses a shortcut, $y$ is closest to $t$), then the
conditional law of the point process in the circle $C(t,r(t,y))$ is
no longer unit rate Poisson. The fact that there are no shortcuts from
$x$ to nodes within this circle biases the probability law and greatly
complicates the analysis. Our approximate greedy algorithm gets around
this problem.

\medskip

Observe that if the message passes through a node $x$, the value of
$r$ immediately after visiting $x$ lies between $\dist(x,t)$ and
$2\dist(x,t)$.

\medskip

We have implicitly assumed that any node can find a local contact closer
to $t$ than itself. We first show that this assumption holds whp if
$c$ is chosen sufficiently large.

\medskip

Fix $c>0$ and $n\in \N$.
For two points $x$ and $y$ in the square $[0,\sqrt{n}]^2$, and a
realisation $\omega$ of the unit rate Poisson process on the square,
define the properties
$$
\cP_n(x,y,\omega) = \{ \exists \ u\in \omega: \dist(u,y) < \dist(x,y)
\quad \mbox{and} \quad \dist(u,x) \le \sqrt{c\log n} \},
$$
and
$$
\cP_n(\omega) = \bigwedge_{(x,y):\dist(x,y) \ge \sqrt{c\log n}}
\cP_n(x,y,\omega).
$$
\begin{lem} \label{lem:good_local_contact}
If $c>0$ is sufficiently large, then
$P(\cP_n(\cdot)) \to 1$ as $n$ tends to infinity.
\end{lem}

In words, with high probability, any two points $x$ and $y$ in the square
$[0,\sqrt{n}]^2$ with $\dist(x,y) > \sqrt{c\log n}$ have the property
that there is a point $u$ of the unit rate Poisson process within
distance $\sqrt{c\log n}$ of $x$ which is closer than $x$ to $y$.
In particular, if $x$ and $y$ are themselves points of the Poisson
process, then $u$ is a local contact of $x$ which is closer to $y$.
The key point to note about the lemma is that it gives a probability
bound which is uniform over all such node pairs.

\begin{proof} Suppose $\dist(x,t)\ge \sqrt{c\log n}$. Consider the
circle $C_1$ of radius $\sqrt{c\log n}$ centred at $x$ and the
circle $C_2$ of radius $\dist(x,t)$ centred at $t$. For any point
$y\neq x$ in their intersection, $\dist(y,t) < \dist(x,t)$.
Moreover, the intersection contains a sector of $C_1$ of angle
$2\pi/3$. Denote this sector $D_1$. Now consider a tessellation of
the square $[0,\sqrt{n}]^2$ by small squares of side $\beta
\sqrt{c\log n}$.
%(A similar construction was used in \cite{xk04}.)
Note that for a sufficiently small geometrical constant $\beta$ that
doesn't depend on $c$ or $n$ ($\beta =1/2$ suffices), the sector $D_1$
fully contains at least one of the smaller squares.
Hence, if every small square contains at least one point of the Poisson
process, then every node at distance greater than $\sqrt{c\log n}$ from
$t$ can find at least one local contact which is closer to $t$.
Number the small squares in some order and let $X_i$ denote the number
of nodes in the $i^{\rm th}$ small square, $i=1,\ldots,n/(\beta^2 c\log n)$.
The number of squares is assumed to be an integer for simplicity.
Clearly, the $X_i$ are iid Poisson random variables with mean
$\beta^2 c\log n$. Hence, by the union bound,
$$
P(\exists \ i: X_i=0) \le \sum_{i=1}^{n/(\beta^2 c\log n)} P(X_i=0)
= \frac{n}{\beta^2 c\log n} e^{-\beta^2 c\log n},
$$
which goes to zero as $n$ tends to infinity, provided that $\beta^2 c>1$.
In particular, $c>4$ suffices since we can take $\beta=1/2$.
\end{proof}

We now state the main result of this section.

\begin{thm}
\label{thm:beta2} Consider the small world random graph described above
with $\alpha =2$, expected node degree $\overline{d}=1$, and $c>0$
sufficiently large, as required by Lemma \ref{lem:good_local_contact}.
Then, the number of hops for message delivery between any pair of nodes
is of order $\log^2 n$ whp.
\end{thm}

\begin{proof} We first evaluate the normalisation constant
$a_n$ by noting that the expected degree, $\overline{d}$, of
a node located at the centre of the square satisfies
$$
\overline{d} \le a_n \int_{\sqrt{c\log n}}^{\sqrt{n/2}}
x^{-2} 2\pi x dx = \pi a_n (\log n-\log \log n - \log (2c)),
$$
and so
\begin{equation} \label{eq:normalisation1}
a_n \ge \frac{1}{\log n},
\end{equation}
for all $n$ sufficiently large, by the assumption that $\overline{d}=1$.

\medskip

Next, we compute the probability of finding a suitable shortcut at
each step of the greedy routing algorithm. We think of the routing
algorithm as proceeding in phases. The value of $r$ is halved at
the end of each phase. The value of $r$ immediately after the
message reaches a node $x$ satisfies the relation $\dist(x,t) \in
(r/2,r]$ at each step of the routing algorithm.
%If $r \le \sqrt{c\log n}$, then $t$ is a local contact of $x$ and
%the message can be routed to $t$ in a single hop. Therefore,
We suppose that $r > k \sqrt{c\log n}$, for some large constant $k$.

\medskip

Denote by $N_A$ the number of nodes in the annulus $A(t,r)$ and observe
that $N_A$ is Poisson with mean $3\pi r^2/16$. The distance from
$x$ to any of these nodes is bounded above by $3r/2$, and so the
probability that a shortcut from $x$ is incident on a particular one
of these nodes is bounded below by $a_n (3r/2)^{-2}$. Thus,
conditional on $N_A$, the probability that $x$ has a shortcut to one
of the $N_A$ nodes in $A(t,r)$ is bounded below by
\begin{equation} \label{eq:good_hop_prob1}
p(r,N_A) = 1 - \Bigl( 1- \frac{4a_n}{9r^2} \Bigr)^{N_A}.
\end{equation}
If $x$ doesn't have such a shortcut, the message is passed via local contacts
which are successively closer to $t$, and hence satisfy the same lower bound
on the probability of a shortcut to $A(t,r)$. Consequently, the number of
local steps $L_x$ until a shortcut is found is bounded above by a
geometric random variable with conditional mean $1/p(r,N_A)$.
Since $N_A \sim \mbox{Pois}(3\pi r^2/16)$, we have by a standard
application of the Chernoff bound that
$$
P(N_A \le \gamma r^2/16) \le \exp \Bigl( -\frac{(3\pi-\gamma)r^2}{16}
+ \frac{\gamma r^2}{16}\log\frac{3\pi}{\gamma} \Bigr),
$$
for any $\gamma<3\pi$.

\medskip

Suppose first that $r \ge k\sqrt{c\log n}$ for some large constant $k$.
Taking $\gamma = 3\pi/2$, we obtain
\begin{equation} \label{eq:node_number_bd}
P \Bigl( N_A \le \frac{3\pi r^2}{32} \Bigr) \le \exp \Bigl(
-\frac{3\pi k^2 c\log n}{32} (1-\log 2) \Bigr).
\end{equation}

Suppose first that $N_A < 3\pi r^2/32$. The number of local hops,
$L_x$, to route the message from $x$ to $A$ is bounded above by the
number of nodes outside $A$, since the distance to $t$ is strictly
decreasing after each hop. Hence,
\begin{equation} \label{eq:localhops_condbd1}
E \Bigl[ L_x \Bigm| N_A < \frac{3\pi r^2}{32} \Bigr] \le
n- \mbox{area}(A) \le n.
\end{equation}

Next, if $N_A \ge 3\pi r^2/32$, then we have by
(\ref{eq:good_hop_prob1}) and (\ref{eq:normalisation1}) that
%\begin{equation} \label{eq:good_hop_prob2}
$$
p(r,N_A) \ge 1 - \exp \Bigl( -\frac{\pi a_n}{24} \Bigr) \ge
1 - \exp \Bigl( -\frac{\pi}{24\log n} \Bigr) \ge \frac{\pi}{48\log n},
$$
%\end{equation}
where the last inequality holds for all $n$ sufficiently large.
Since the number of hops to reach $A$ is bounded above by a geometric
random variable with mean $1/p(r,N_A)$, we have
\begin{equation} \label{eq:localhops_condbd2}
E \Bigl[ L_x \Bigm| N_A \ge \frac{3\pi r^2}{32} \Bigr] \le
\frac{48}{\pi}\log n.
\end{equation}
Finally, we obtain from (\ref{eq:node_number_bd}),
(\ref{eq:localhops_condbd1}) and (\ref{eq:localhops_condbd2}) that
$$
E[L_x] \le n\exp \Bigl( -\frac{3\pi k^2 c(1-\log 2)}{32} \log n \Bigr)
+ \frac{48}{\pi}\log n.
$$
The first term in the sum above can be made arbitrarily small by
choosing $k$ large enough, so $E[L_x]=O(\log n)$. It can also be seen
from the arguments above that $L_x=O(\log n)$ whp. In other words,
while $r\ge k\sqrt{c\log n}$, the number of hops during each phase
is of order $\log n$. Moreover, the number of such phases is of order
$\log n$ since the initial value of $r$ is at most $\sqrt{2n}$, and
$r$ halves at the end of each phase.

\medskip

Hence, the total number of hops until $r < k\sqrt{c\log n}$ is of
order $\log^2 n$. Once the message reaches a node $x$ with $\dist(x,t)
< k\sqrt{c\log n}$, the number of additional hops to reach $t$ is bounded
above by the total number of nodes in the circle $C(t,k\sqrt{c\log n})$.
By using the Chernoff bound for a Poisson random variable, it can be shown
that this number is of order $\log n$ whp. This completes the
proof of the theorem.
\end{proof}

\section{Impossibility of efficient routing when $\alpha \neq 2$}

We now show that if $\alpha<2$, then no decentralised algorithm can
route between arbitrary source-destination pairs in time which is
polylogarithmic in $n$. In fact, the number of routing hops is
polynomial in $n$ with some fractional power that depends on $\alpha$.

\medskip

We now make precise what we mean by a decentralised routing
algorithm. As specified earlier, each node knows the locations of
all its local contacts with distance $\sqrt{c\log n}$ and of all
its shortcut neighbours, as well as other nodes (if any) from which
shortcuts are incident to it. A routing algorithm
specifies a (possibly random) sequence of nodes $s=x_0, x_1,
\ldots, x_k = t, x_{k+1}=t, \ldots$, where the only requirement is
that each node $x_i$ be chosen from among the local or shortcut
contacts of nodes $\{ x_0,\ldots,x_{i-1} \}$. (This is the same
definition as used by Kleinberg \cite{Klei00}).

\begin{thm}
\label{thm:betaless2} Consider the small world random graph
described above with $\alpha < 2$, and arbitrarily large constants $c$
and $\overline{d}$. Suppose the source $s$ and destination $t$ are
chosen uniformly at random from the node set. Then,
%there is a $\gamma>0$ such that
the number of hops for message delivery in any decentralised
algorithm exceeds $n^{\gamma}$ whp, for any $\gamma < (2-\alpha)/6$.
\end{thm}

It is not important that the source and destination be chosen uniformly
but only that the distance between them be of order $n^a$ whp for some
$a>0$.

\begin{proof} We first evaluate the normalisation constant
$a_n$ by noting that the expected degree satisfies
$$
\overline{d} \ge a_n \int_{\sqrt{c\log n}}^{\sqrt{n}/2}
x^{-\alpha} 2\pi x dx = \frac{2 \pi a_n}{2-\alpha} \Bigl(
\frac{n^{(2-\alpha)/2}}{2^{2-\alpha}} - (c\log n)^{(2-\alpha)/2} \Bigr),
$$
which, on simplification, yields that
\begin{equation} \label{eq:normalisation2}
a_n \le \frac{4\overline{d}}{n^{(2-\alpha)/2}},
\end{equation}
for all $n$ sufficiently large. Note that $a_n$ is an upper bound
on the probability that there is a shortcut between any pair of nodes.

\medskip

Suppose that the source $s$ and destination $t$ are chosen uniformly from
all nodes on $[0,\sqrt{n}]^2$. Fix $\delta \in (\gamma,1/2)$ and define
$C_{\delta} = C(t,n^{\delta})$ to be the circle of radius $n^{\delta}$
centred at $t$. It is clear that, for any $\epsilon>0$, the distance
$\dist(s,C_{\delta})$ from $s$ to the circle $C_{\delta}$ is bigger than
$n^{(1/2)-\epsilon}$ whp. Suppose now that this inequality holds, but
that there is a routing algorithm which can route from $s$ to $t$ in
fewer than $n^{\gamma}$ hops. Denote by $s=x_0, x_1, \ldots, x_m=t$,
the sequence of nodes visited by the routing algorithm, with
$m\le n^{\gamma}$. We claim that there must be a shortcut from at least
one of the nodes $x_0,x_1,\ldots, x_{m-1}$ to the set $C_{\delta}$.
Indeed, if there is no such shortcut, then $t$ must be reached starting
from some node outside $C_{\delta}$ and using only local links. Since
the length of each local link is at most $\sqrt{c\log n}$ and the number
of hops is at most $n^{\gamma}$, the total distance traversed by local
hops is strictly smaller than $n^{\delta}$ (for large enough $n$, by
the assumption that $\delta>\gamma$), which yields a contradiction.
We now estimate the probability that there is a shortcut from one of
the nodes $x_0,\ldots,x_{m-1}$ to the set $C_{\delta}$.

%Hence, for any two nodes $u$ and $v$
%which are more than $\sqrt{c\log n}$ apart,
%the probability that there is no shortcut between them,
%i.e., none
%of the $2\overline{d}$ shortcuts from these nodes are incident on the
%other,
%is bounded below by
%$$
%(1-a_n)^{2\overline{d}} \ge 1-2\overline{d}a_n.
%$$
%Hence, for all nodes $u$ and $v$,
%\begin{equation} \label{eq:shortcut_bd1}
%P(\exists \mbox{ shortcut between $u$ and $v$}) \le \frac{8}{n^{(2-\alpha)/2}}.
%\end{equation}

%Fix $\delta \in (0,1/2)$ and define $C_{\delta} = C(t,n^{\delta})$ to be
%the circle of radius $n^{\delta}$ centred at $t$.

\medskip

The number of nodes in the circle $C_{\delta}$, denoted $N_C$,
is Poisson with mean $\pi n^{2\delta}$, so $N_C < 4n^{2\delta}$ whp.
Now, by (\ref{eq:normalisation2}) and the union bound,
$$
P(\exists \mbox{ shortcut between $u$ and $C_{\delta}$}|N_C <
4n^{2\delta}) \le 16 \ \overline{d} \ n^{(4\delta+\alpha-2)/2},
$$
for any node $u$. Applying this bound repeatedly for each of the nodes
$x_0, x_1, \ldots, x_{m-1}$ generated by the routing algorithm, we get,
\begin{equation} \label{eq:shortcut_bd2}
P(\exists \mbox{ shortcut to $C_{\delta}$ within $n^{\gamma}$
hops)} |N_C < 4n^{2\delta}) \le 16 \ \overline{d} \ n^{(2\gamma +
4\delta+\alpha-2)/2}.
\end{equation}
%where $\pi_k := \{ x_0,x_1,\ldots,x_{k_1} \}$. Take $k=n^{\gamma}$;
%since $\gamma < (2-\alpha)/2$ by assumption, we can choose $\delta>0$
Now $\gamma < (2-\alpha)/6$ by assumption, and $\delta>\gamma$ can be
chosen arbitrarily. In particular, we can choose $\delta$ so that
$2\gamma + 4\delta + \alpha -2$ is strictly negative, in which case the
conditional probability of a shortcut to $C_{\delta}$ goes to zero
as $n\to \infty$. Since $P(N_C \ge 4n^{2\delta})$ also goes to
zero, we conclude that the probability of finding an $s-t$ route
with fewer than $n^{\gamma}$ hops also goes to zero. This
concludes the proof of the theorem.
\end{proof}

{\bf Remarks:} The theorem continues to hold if we assume 1-step lookahead.
By this, we mean that when a node decides where to send the message at the
next step, it can not only use the locations of all its local and shortcut
contacts, but also the locations of their contacts. All this means is that
after visiting $n^{\gamma}$ nodes, the algorithm has knowledge about
$O(n^{\gamma} \log n)$ nodes. If none of these nodes has a shortcut into
the set $C_{\delta}$, which is the case whp, then the arguments above still
apply. The same is true for $k$-step lookahead, for any constant $k$.

\begin{thm}
\label{thm:betamore2} Consider the small world random graph
described above with $\alpha > 2$, and arbitrarily large constants $c$
and $\overline{d}$. Suppose the source $s$ and destination $t$ are
chosen uniformly at random from the node set. Then,
%there is a $\gamma>0$ such that
the number of hops for message delivery in any decentralised
algorithm exceeds $n^{\gamma}$ whp, for any $\gamma <
(\alpha-2)/(2(\alpha-1))$.
\end{thm}

\begin{proof} For a node $u$, the probability that a randomly
generated shortcut has length bigger than $r$ is bounded above by
$$
\frac{ \int_r^{\infty} x^{-\alpha} 2\pi xdx }{ \int_{\sqrt{c\log
n}}^{\sqrt{n}/2} x^{-\alpha} 2\pi xdx} \le \mbox{const. }r^{2-\alpha} (\log
n)^{(\alpha-2)/2},
$$
for all $n$ sufficiently large. Since there are $2\overline{d}$ shortcuts
per node on average, the probability that two nodes $u$ and $v$ separated
by distance $r$ or more possess a shortcut between them is bounded above
by the same function, but with the constant suitably modified.

\medskip

Now, for randomly chosen nodes $s$ and $t$, $\dist(s,t) > n^{(1/2)-\epsilon}$
whp, for any $\epsilon>0$. Hence, there can be a path of length $n^{\gamma}$
hops between $s$ and $t$ only if at least one of the hops is a shortcut
of length $n^{(1/2)-\epsilon-\gamma}$ or more. By the above and the union
bound, the probability of there being such a shortcut is bounded above by
$$
\mbox{const. } n^{\gamma} \Bigl( n^{(1/2)-\epsilon-\gamma}
\Bigr)^{2-\alpha} (\log n)^{(\alpha-2)/2}.
$$
The exponent of $n$ in the above expression is
$$
\frac{2-\alpha}{2}(1-2\epsilon) + \gamma(\alpha-1).
$$
The exponent above is negative for sufficiently small $\epsilon>0$
provided $\gamma < (\alpha-2)/(2(\alpha-1))$. In other words, if this
inequality is satisfied, then the probability of finding a route
with fewer than $n^{\gamma}$ hops goes to zero as $n\to \infty$.
This establishes the claim of the theorem.
\end{proof}

%\newpage


\begin{thebibliography}{10}

\bibitem{BR}
A.~D.~Barbour and G.~Reinert, ``Small worlds", \emph{Random Structures
and Algorithms}, 19: 54--74, 2001.

\bibitem{BB}
I.~Benjamini and N.~Berger, ``The diameter of long-range percolation
clusters on finite cycles", \emph{Random Structures and Algorithms},
19: 102--111, 2001.

\bibitem{BKPS}
I.~Benjamini, H.~Kesten. Y.~Peres and O.~Schramm, ``Geometry of the
uniform spanning forest: phase transitions in dimensions 4,8,12,\ldots",
\emph{Annals of Mathematics} 160: 465--491, 2004.

\bibitem{CGS}
Don Coppersmith, David Gamarnik and Maxim Sviridenko, ``The diameter
of a long-range percolation graph", {\em Random Structures and Algorithms},
21: 1--13, 2002.

\bibitem{FrMee04}
M.~Franceschetti and R.~Meester,
\newblock ``Navigation in small world networks, a scale-free continuum
approach".
\newblock {\em Preprint}, 2004.

\bibitem{Klei00}
J.~M.~Kleinberg,
\newblock      ``The small-world phenomenon: an algorithmic perspective",
\newblock {\em Proc. 32nd Annual ACM Symposium on the Theory of Computing
(STOC)}: 163--170, 2000.

\bibitem{Mil67}
S.~Milgram,
\newblock ``The small world problem",
\newblock {\em Psychology {T}oday}, 2: 60--67, 1967.

\bibitem{penrose}
Mathew Penrose,
\newblock \emph{Random Geometric Graphs},
\newblock Oxford University Press, 2003.

\bibitem{SM}
G.~Sharma and R.~R.~Mazumdar, ``Hybrid sensor networks: A small world",
\emph{Proc. ACM MobiHoc}, 2005.

\bibitem{WaSt98}
D.~J.~Watts and S.~H.~Strogatz,
\newblock   ``Collective dynamics of small world networks",
\newblock {\em Nature}, 393: 440--442, 1967.

\end{thebibliography}
\end{document}